\documentclass[12pt]{article}
\usepackage{epsfig} 
 
\begin{document}
 
\def\bea{\begin{eqnarray}} \def\eea{\end{eqnarray}} \def\Pl{{\cal P}(l)} 

\title{\Large \bf Separations inside a cube \\ } 
\author{A.F.F. Teixeira \thanks{teixeira@cbpf.br} \\ 
{\small Centro Brasileiro de Pesquisas F\'{\i}sicas } \\ 
{\small 22290-180 Rio de Janeiro-RJ, Brazil} \\ } 
\date{\today} 
\maketitle

\begin{abstract} 
Two points are randomly selected inside a three-dimensional euclidian cube. The value $l$ of their separation lies somewhere between zero and the length of a diagonal of the cube. The probability density ${\cal P}(l)$ of the separation is constructed analytically. Also a Monte Carlo computer simulation is performed, showing good agreement with the formulas obtained. 
\end{abstract}

\section{Introduction} 
An important problem in geometry and statistics is: given a convex compact space endowed with a metric, and randomly choosing two points in the space, find the probability density $\Pl$ that these points have a specified separation $l$. The study of this problem has a long history \cite{Beres}, and recently gained considerable impetus from researchers in cosmic crystallography \cite{FarrarMelott}-\cite{Janna}. 

In a recent paper the functions $\Pl$ corresponding to $2D$ disks and rectangles were obtained \cite{Memb}. The methodology introduced in that work is here extended to a $3D$ euclidian cube.

\section{Preliminaries} 
An euclidian cube with side $a$ is assumed, occupying the location $0<x,y,z<a$ in a cartesian frame. Randomly choosing two points $A$ and $B$ in the cube, we want the probability $\Pl dl$ that the separation between the points lie between $l$ and $l+dl$. The probability density $\Pl$ has to satisfy the normalization condition 
\bea 											\label{a1}
\int_0^{\sqrt{3}a}\Pl dl =1 . 
\eea 
The calculation can be shortened if one considers the symmetries of the cube. Really, if the points $A$ and $B$ have been chosen, imagine the oriented segment $A'B'$ parallel to $AB$, with the tip $A'$ coinciding with the origin $O$. The other tip $B'$ then lies inside a larger cube, with side $2a$. Since the probability density $\Pl$ clearly does not depend on which octant of the large cube contains $B'$, there is no loose in generality in restricting the calculation to the cases where $B'$ is in the octant $0<x,y,z<a$. 

With this assumption, the point $B'$ has cartesian coordinates 
\bea \nonumber 
B'=(l\cos\theta\cos\phi,\hspace*{1mm} l\cos\theta\sin\phi,\hspace*{1mm} l\sin\theta), 
\eea
where both angles $\theta, \phi$ are bound to the interval $[0, \pi/2]$; here $\phi$ is the azimuthal angle, while $\theta$ is the polar angle measured from the $z=0$ plane. The corresponding tip $B$ in the original segment must lie inside a parallelepiped with sides (see figure 1) 
\bea											\label{a2} 
l_x:=a-l\cos\theta\cos\phi,\hspace*{3mm} l_y:=a-l\cos\theta\sin\phi,\hspace*{3mm} l_z:=a-l\sin\theta. 
\eea 

\vspace*{3mm}
\centerline{\epsfig{file=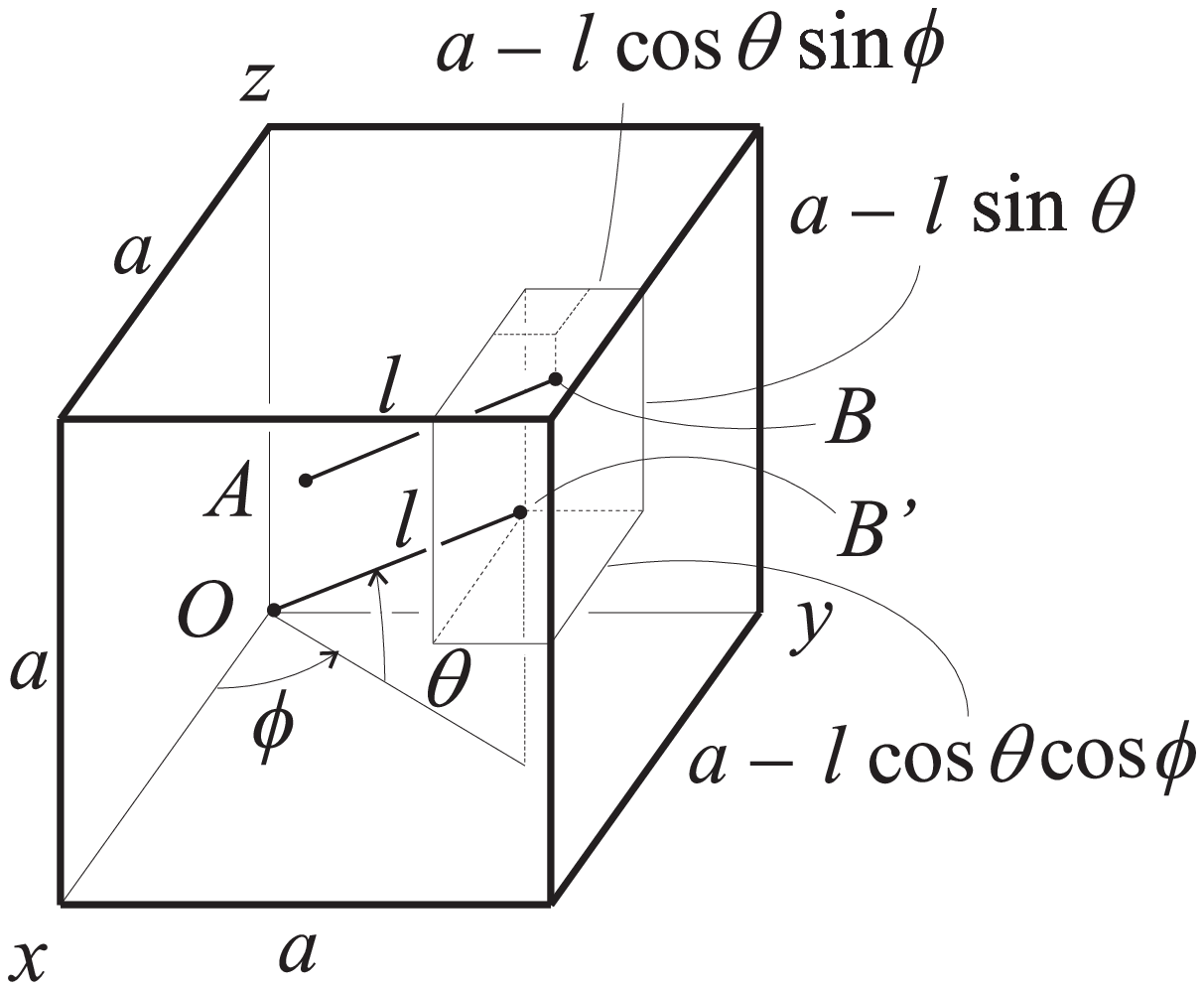,width=6cm,height=5cm}} 
\vspace*{3mm} 

\noindent {\small {\bf Figure 1} The endpoint $B$ of the segment $AB$ must lie inside the parallelepiped with a corner at $B'(l,\phi,\theta$).} 
\vspace*{5mm}

The probability ${\cal P}(l,\theta,\phi) dl\, d\theta\, d\phi$ that the segment $AB$ has length between $l$ and $l+dl$, azimuth between $\phi$ and $\phi+d\phi$, and polar angle between $\theta$ and $\theta+d\theta$ is then 
\bea 											\label{a3} 
{\cal P}(l,\theta,\phi) dl\, d\theta\, d\phi = k\,l_x\,l_y\,l_z\, l^2\cos\theta\, dl\, d\theta\, d\phi , 
\eea 
where $k$ is a constant and where the assumption $0<\phi,\theta<\pi/2$ stands. Performing the angular integrations we shall obtain 
\bea 											\label{a4} 
\Pl = \int\!\!\int{\cal P}(l,\theta,\phi) d\theta\,\, d\phi, 
\eea 
and we finally fix $k$ using the condition (\ref{a1}).
 
To calculate $\Pl$, three cases need be separately considered, depending on the value of $l$ relative to $a$: namely the cases $0< l< a,$  $a< l<\sqrt{2}a,$ and $\sqrt{2}a< l<\sqrt{3}a$.

\section{The case $0< l< a$}  
As is seen in the figure 2, in this case we effectively have $\phi_{min}=\theta_{min}=0,$ and 
$\phi_{max}=\theta_{max}=\pi/2.$ Then 
\bea 											\label{a5} 
{\cal P}(l< a)=k\,l^2\int_0^{\pi/2}l_z \cos\theta\, d\theta\int_0^{\pi/2}l_x\,l_y\, d\phi\\ 
\label{a6}  =\frac{k\,l^2}{8}[4\pi a^3-6\pi a^2l+8al^2-l^3],  \eea 
where $k=8/a^6$ as will be fixed later on. 

\vspace*{3mm}
\centerline{\epsfig{file=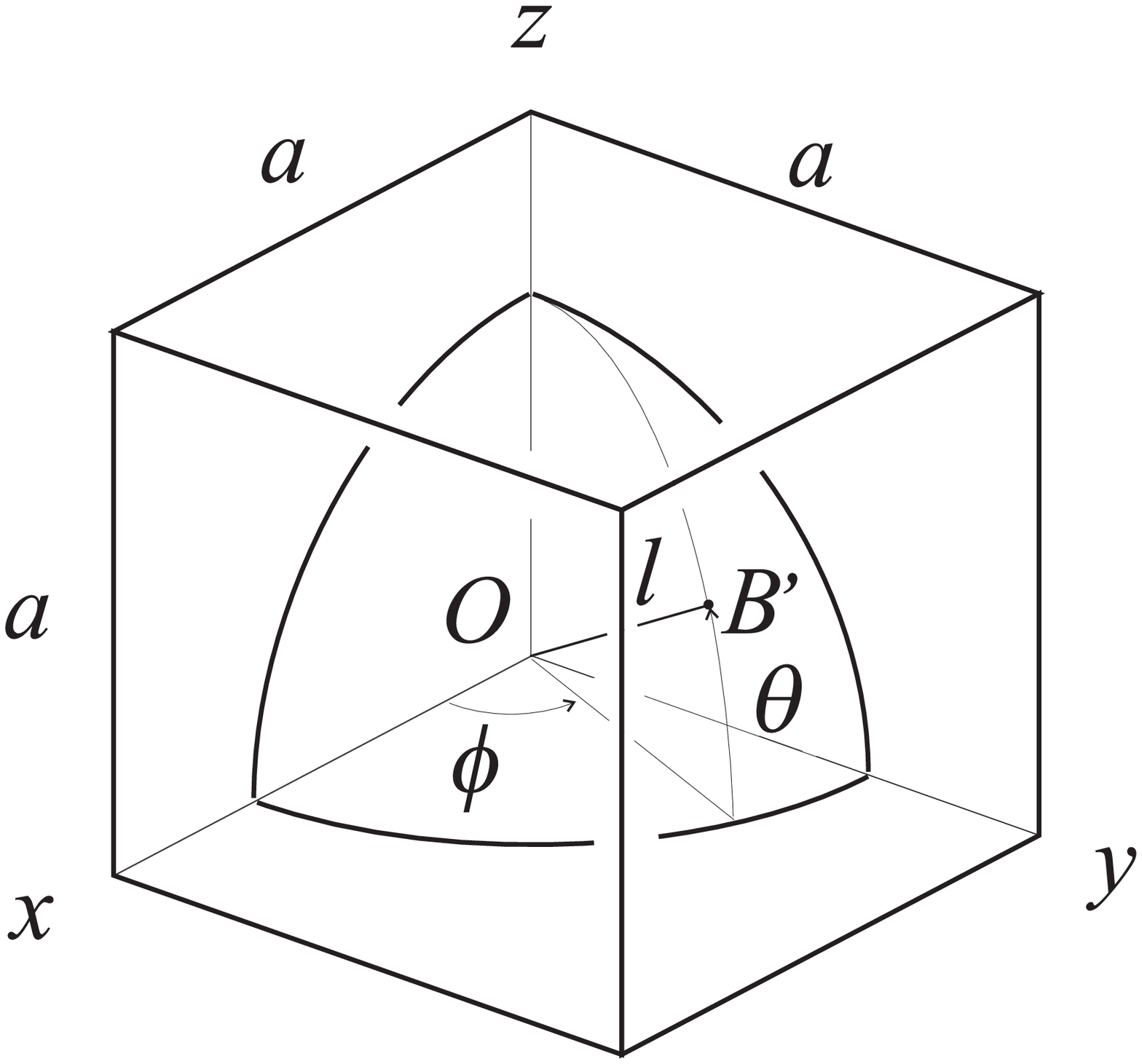,width=5cm,height=5cm}} 
\vspace*{3mm} 

\noindent {\small {\bf Figure 2} The triangular intersection of the $2D$ sphere with radius $l$ and centre $O$ with the $3D$ cube having side $a>l$. }  
\vspace*{5mm}

\section{The case $a< l< \sqrt{2}a$}  
In this case the intersection of the $2D$ sphere (with radius $l$) with the $3D$ cube (with side $a$) is an hexagonal surface as in figure 3. 

\vspace*{3mm} 
\centerline{\epsfig{file=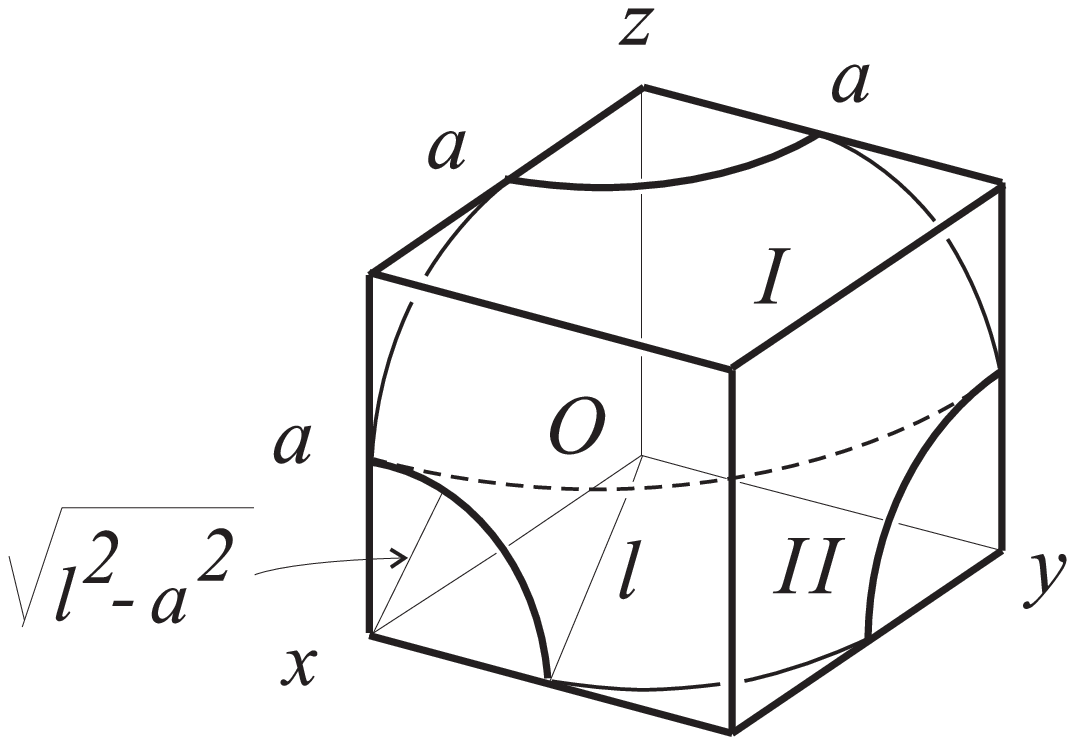,width=6cm,height=45mm}} 
\vspace*{3mm} 

\noindent{\small {\bf Figure 3} The hexagonal intersection of a $2D$ sphere with radius $l$ and centre $O$ with a $3D$ cube with a vertex in $O$ and having side $a$ such that $a<l<\sqrt{2}a$.}  
\vspace*{5mm}

\noindent We note that the arcs of circle drawn on the faces $x,y,z=0$ have radius $l$, while those drawn on the faces $x,y,z=a$ have radius $\sqrt{l^2-a^2}$. 

For convenience of integration we divide the intersection into two regions. In region $I$ we have $\cos^{-1}(a/l)<\theta<\sin^{-1}(a/l)$ and $0<\phi<\pi/2$. 

In region $II$ we have $\theta_{min}=0$ and $\theta_{max}=\cos^{-1}(a/l)$. To have $\phi_{min}(\theta)$ we note that the circle drawn on the face $x=a$ satisfies the equation $\cos\phi\,\cos\theta=a/l$, so 
\bea 											\label{a7}
\phi_{min}(\theta)=\cos^{-1}\left(\frac{a}{l\cos\theta}\right)=:\phi_1(\theta). 
\eea
\noindent On the other hand, the circle drawn on the face $y=a$ satisfies $\sin\phi\,\cos\theta=a/l$, so we have 
\bea 											\label{a8}
\phi_{max}(\theta)=\sin^{-1}\left(\frac{a}{l\cos\theta}\right)=:\phi_2(\theta). 
\eea 
We then find  
\bea 									
{\cal P}(a<l<\sqrt{2}a) \nonumber  
\eea 
\bea
=k\,l^2\left[\int_{\cos^{-1}(a/l)}^{\sin^{-1}(a/l)}\cos\theta d\theta \int_0^{\pi/2}d\phi + \int_0^{\cos^{-1}(a/l)}\cos\theta d\theta\int_{\phi_1(\theta)}^{\phi_1(\theta)}d\phi\right]l_xl_yl_z   \label{a9} \\ 
=\frac{k\,l}{8}\Big[ 2l^4+6a^2l^2-a^4-2\pi a^3(4l-3a)-8a(2l^2+a^2)\sqrt{l^2-a^2} \nonumber \\ +24a^2l^2\cos^{-1}(a/l)\Big], \label{a10}
\eea 
where $k=8/a^6$ as will be fixed later on.

\section{The case $\sqrt{2}a<l<\sqrt{3}a$} 
In this case the $2D$ sphere with radius $l$ intersects the $3D$ cube with side $a$ in the triangular surface shown in figure 4. 

\vspace*{3mm} 
\centerline{\epsfig{file=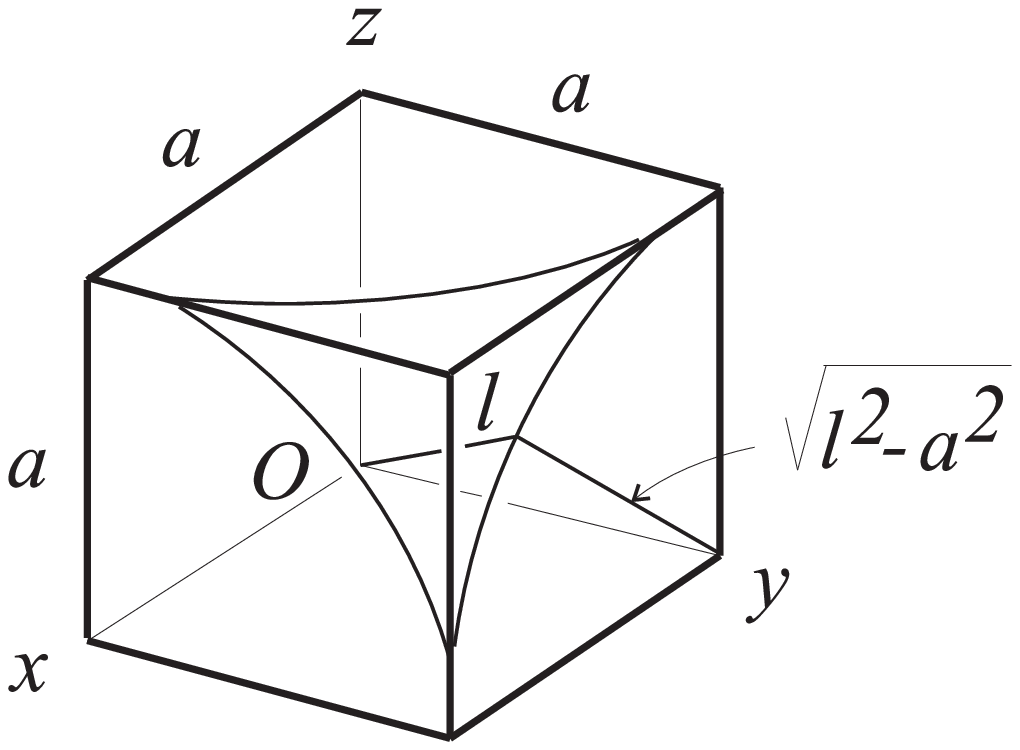,width=6cm,height=45mm}}  
\vspace*{3mm} 

\noindent {\small {\bf Figure 4} The triangular intersection of a $2D$ sphere with radius $l$ and centre $O$ with a $3D$ cube with a vertex in $O$ and having side $a$ such that $\sqrt{2}a<l<\sqrt{3}a$. }  
\vspace*{5mm}

\noindent As before, the circles drawn on the faces $x,y,z=a$ have radius $\sqrt{l^2-a^2}$. The azimuthal integration is performed between $\phi_1(\theta)$ and 
$\phi_2(\theta)$ as in the region $II$ of the preceding case, and again $\sin\theta_{max}=a/l$; but now $\cos\theta_{min}=\sqrt{2}a/l$. We then find 
\bea											 
{\cal P}(\sqrt{2}a<l<\sqrt{3}a)=  
k\,l^2\int_{\cos^{-1}(\sqrt{2}a/l)}^{\sin^{-1}(a/l)}l_z\cos\theta d\theta\int_{\phi_1(\theta)}^{\phi_2(\theta)}l_xl_yd\phi             \label{a11}\\ 
=\frac{k\,l}{8}\Big[ 8a(l^2+a^2)\sqrt{l^2-2a^2}-(l^2+a^2)(l^2+5a^2) +2\pi a^2(3l^2-4al+3a^2) \nonumber \\ +24a^3l\sec^{-1}(l^2/a^2-1)-24a^2(l^2+a^2)\sec^{-1}\sqrt{l^2/a^2-1}\Big] ,                                                                      \label{a12}
\eea 
\noindent where $k=8/a^6$. 

This value for the constant $k$ derives from the normalization condition (\ref{a5}), namely, 
\bea										         \label{a13} 
\int_0^a\!\!{\cal P}(l<a)dl+\int_a^{\sqrt{2}a}\!\!{\cal P}(a<l<\sqrt{2}a)dl+ \int_{\sqrt{2}a}^{\sqrt{3}a}\!\!{\cal P}(\sqrt{2}a<l<\sqrt{3}a)dl=1 .
\eea

\section{Graphs of $\Pl$} 
In figure 5 we present a graph of the dimensionless function $a\Pl$ against the dimensionless variable $l/a$. 

\vspace*{3mm} 
\centerline{\epsfig{file=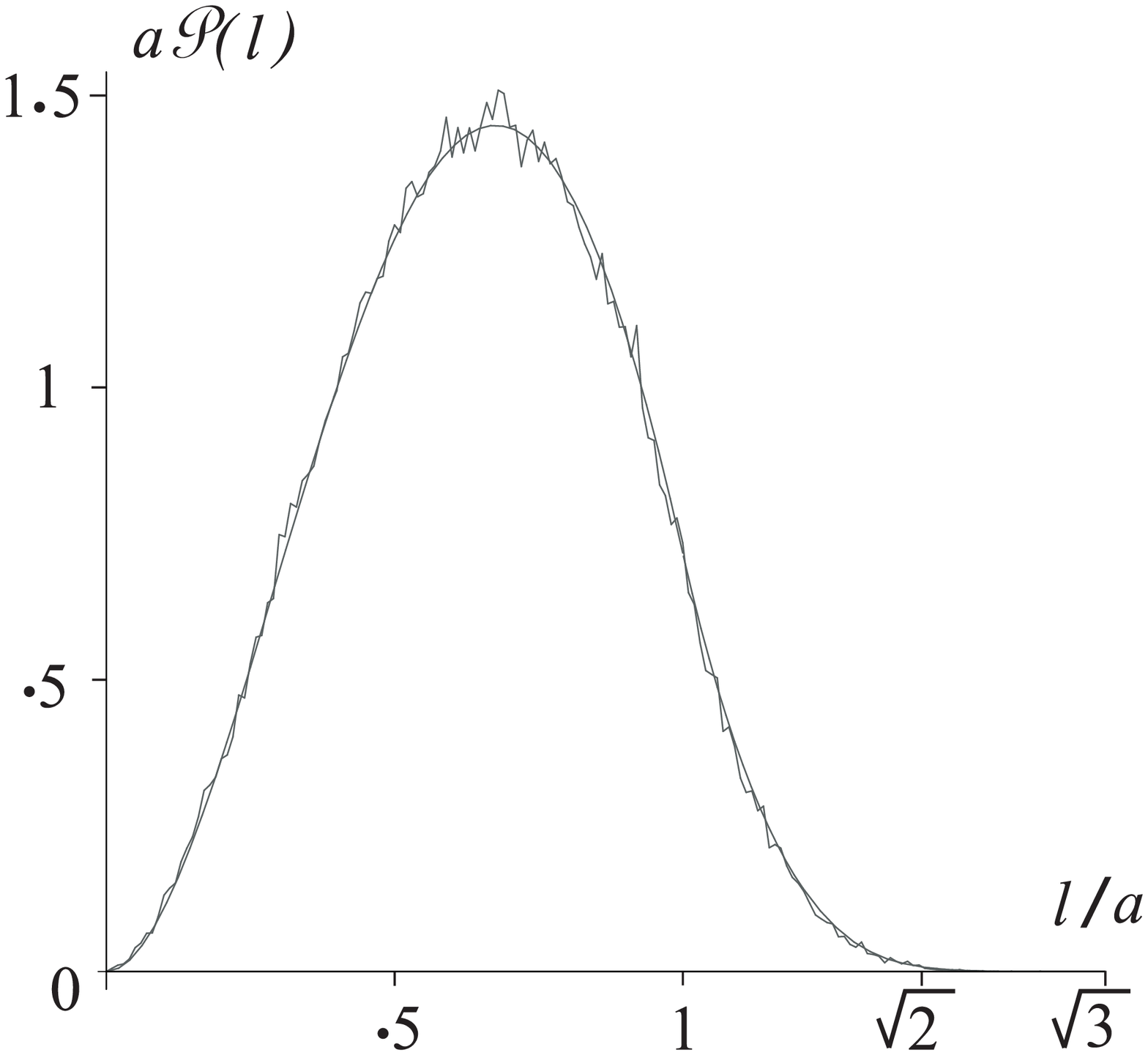,width=6cm,height=55mm}} 
\vspace*{3mm} 

\noindent {\small {\bf Figure 5} The probability density $\Pl$ of separation $l$ of pairs of randomly distributed points inside a cube with side $a$. The irregular curve is the output of a corresponding computer simulation. }  
\vspace*{5mm}

\noindent We note that the function and its first derivative are continuous in the whole interval $0<l<\sqrt{3}a$. Nevertheless the second derivative is discontinuous at $l=a$, as discussed in the next section. In the figure a normalized histogram corresponding to 150,000 separations between pairs of points randomly selected in the cube is superimposed, for comparison; the agreement of the two curves evinces the correctness of the calculation.

\section{Comments} 
The integration to find $\Pl$ in eqs. (\ref{a5})-(\ref{a6}) is almost trivial; however, not the same can be said about the two other cases, namely in going from (\ref{a9}) to (\ref{a10}) and from (\ref{a11}) to (\ref{a12}). A computer assistance appears paramount in these two cases, to confirm every short step in the calculation and simplification of expressions. 

Similarly as in \cite{Memb}, the probability density $\Pl$ and its first derivative are continuous throughout the entire range $0<l<\sqrt{3}a$. But the second derivative shows a finite discontinuity at $l=a$, although it is continuous at $l=\sqrt{2}a$. 

A remarkable feature of $\Pl$ is its behaviour for large values of $l$; really, near $l=\sqrt{3}a$ we find 
\bea 											\label{a14}
a\Pl=\frac{9}{5}(\sqrt{3}-l/a)^5+O\big((\sqrt{3}-l/a)^6\big),  
\eea 
so $\Pl$ is essentially a fifth power of $\sqrt{3}-l/a$. We find that 91\% of the separations lie in the range $l\in(0, a)$, 9\% lie in the interval $l\in(a, \sqrt{2}a)$, and only 0.04\% have  $l>\sqrt{2}a$.


\begin{thebibliography}{30} 
\bibitem{Beres} Krzysztof Bere\'{s}, {\it Distance distribution}, {\cal Zeszyty Naukowe 
        Universytetu Jagiello\'{n}skiego - Acta Cosmologica - Z. 5} (1976) 7-27 
\bibitem{FarrarMelott} Kelly A. Farrar \& Adrian L. Melott, {\it Gravity in twisted                         space}, {\cal Computers in Physics}, Mar/Apr 1990, 185-189
\bibitem{LeLaLu} Roland Lehoucq, M. Lachi\`{e}ze-Rey \& Jean-Pierre Luminet, {\it Cosmic 
        crystallography}, gr-qc/9604050 
\bibitem{FagGaus1} Helio V. Fagundes \& Evelise Gausmann, {\it On closed Einstein-de Sitter 
        universes}, astro-ph/9704259 
\bibitem{LeLuUz} Roland Lehoucq, Jean-Pierre Luminet \& Jean-Philippe Uzan, {\it Topological 
        lens effects in universes with non-euclidian compact spatial sections}, 
        astro-ph/9811107 
\bibitem{GTRB1} Germ\'{a}n I. Gomero, Antonio F.F. Teixeira, Marcelo J. Rebou\c{c}as 
        \& Armando  Bernui, {\it Spikes in cosmic crystallography}, gr-qc/9811038 
\bibitem{FagGaus2} Helio V. Fagundes \& Evelise Gausmann, {\it Cosmic crystallography in 
        compact hyperbolic universes}, astro-ph/9811368 
\bibitem{LuRouk} Jean-Pierre Luminet \& Boudewijn F. Roukema, {\it Topology of the universe: 
         theory and observation}, astro-ph/9901364 
\bibitem{UzLeLu} Jean-Philippe Uzan, Roland Lehoucq \& Jean-Pierre Luminet, {\it A new 
         crystallographic method for detecting space topology}, astro-ph/9903155
\bibitem{BT1} Armando Bernui \& Antonio F.F. Teixeira, {\it Cosmic crystallography: 
        three multipurpose functions}, astro-ph/9904180
\bibitem{GRT1} Germ\'{a}n I. Gomero, Marcelo J. Rebou\c{c}as \& Antonio F.F. Teixeira, 
       {\it Spikes in cosmic crystallography II: topological signature of compact flat 
         universes}, gr-qc/9909078
\bibitem{GRT2} Germ\'{a}n I. Gomero, Marcelo J. Rebou\c{c}as \& Antonio F.F. Teixeira, 
        {\it A topological signature in cosmic topology}, gr-qc/9911049
\bibitem{BT2} Armando Bernui \& Antonio F.F. Teixeira, {\it Cosmic crystallography: 
         the euclidian isometries}, gr-qc/0003063 
\bibitem{T1} Antonio F.F. Teixeira, {\it Cosmic crystallography in a circle}, 
         gr-qc/0005052  
\bibitem{LeUzLu2} Roland Lehoucq, Jean-Pierre Luminet \& Jean-Philippe Uzan, 
        {\it Limits of crystallographic methods for detecting space topology}, 
        astro-ph/0005515
\bibitem{T2} Antonio F.F. Teixeira, {\it Cosmic crystallography: the hyperbolic                 isometries}, gr-qc/0010107 
\bibitem{GRT3} Germ\'{a}n I. Gomero, Marcelo J. Rebou\c{c}as \& Antonio F.F. Teixeira,  
        {\it Signature for the shape of the universe}, gr-qc/0105048
\bibitem{EvLeLuUzWe} Evelise Gausmann, Roland Lehoucq, Jean-Pierre Luminet, Jean-Philippe 
         Uzan \& Jeffrey Weeks, {\it Topological lensing in spherical spaces}, 
         gr-qc/0106033
\bibitem{Janna} Janna Levin, {\it Topology and the cosmic microwave background}, 
         gr-qc/0108043
\bibitem{Memb} Antonio F.F. Teixeira, {\it Distances in plane membranes}, physics/0111108

\end{thebibliography}
\end{document}